\documentstyle[11pt,twoside]{article}
\newtheorem{theorem}{\noindent\bf Theorem}[section]
\newtheorem{proposition}[theorem]{\noindent\bf Proposition}
\newtheorem{lemma}[theorem]{\noindent\bf Lemma}
\newtheorem{corollary}[theorem]{\noindent\bf Corollary}

\newtheorem{example}[theorem]{\noindent\bf Example}

\newcommand{\A}{{\cal A}}
\newcommand{\h}{{\cal H}}
\newcommand{\m}{\verb"m"}
\newcommand{\M}{\verb"M"}
\newcommand{\e}{\hfill\blacksquare}
\input{amssym}
\begin{document}
\date{}
\title{{\Large\bf Some properties of generalized and approximately\\ dual
 frames in Hilbert spaces}}

\author{{\normalsize\sc H. Javanshiri}}

\maketitle
\normalsize

\begin{abstract}
In the present paper, some sufficient and necessary conditions for
two frames $\Phi=(\varphi_n)_n$ and $\Psi=(\psi_n)_n$ under which
they are approximately or generalized dual frames are determined
depending on the properties of their analysis and synthesis
operators. We also give a new characterization for approximately
dual frames associated with a given frame and given operator by
using of bounded operators. Among other things, we prove that if
two frames $\Phi=(\varphi_n)_n$ and $\Psi=(\psi_n)_n$ are close to
each other, then we can find approximately dual frames
$\Phi^{ad}=(\varphi^{ad}_n)_n$ and $\Psi^{ad}=(\psi^{ad}_n)_n$ of
them which are close to each other and $T_\Phi
U_{\Phi^{ad}}=T_\Psi U_{\Psi^{ad}}$, where $T_\Phi$ and $T_\Psi$
(resp. $U_{\Phi^{ad}}$ and $U_{\Psi^{ad}}$) are the analysis
operators (resp. synthesis operators) of the frames $\Phi$ and
$\Psi$ (resp. $\Phi^{ad}$ and $\Psi^{ad}$), respectively. We then
give some consequences on generalized dual frames. Finally, we
apply these results to find some construction results for
approximately dual frames for a given Gabor frame.\\

{\bf Mathematics Subject Classification}: Primary: 42C15;
Secondary: 47A58.

{\bf Key words}: Frame, Gabor frame, dual frame, approximately
dual frame, generalized

 dual frame.

\end{abstract}


\section{\large\bf Introduction and some prerequisites }

Throughout this paper, we denote by $\cal H$ a separable Hilbert
space with the inner product ``$\big<\cdot,\cdot\big>"$ and for
basic notations and terminologies on the theory of frames for
$\cal H$, we shall follow \cite{c}. Recall that a sequence
$\Phi=(\varphi_n)_n\subseteq{\cal H}$ is a {\it frame} for $\cal
H$, if there exist constants $\m_\Phi, \M_\Phi>0$ such that for
all $f\in\h$
\begin{equation}\label{01}
\m_\Phi\|f\|^2\leq\sum_{n=1}^\infty|\big<f,\varphi_n\big>|^2\leq\M_\Phi\|f\|^2,
\end{equation}
where $\m_\Phi, \M_\Phi$ are called frame bounds. Moreover, the
sequence $\Phi=(\varphi_n)_n$ is called a {\it Bessel sequence}
for $\h$, if only the second inequality of (\ref{01}) holds.
 Let also, the space $\ell^2$
is defined as usual and we denote its canonical orthonormal basis
by $\Delta=(\delta_n)_n$. Moreover, the notation $B({\mathcal H})$
(resp. $B({\mathcal H},\ell^2)$) is used to denote the collection
of all bounded linear operators from $\mathcal H$ into $\mathcal
H$ (resp. $\ell^2$).

It is well-known that the Bessel sequence $\Phi=(\varphi_n)_n$
gives three operators which plays a crucial role in the theory of
frames; In what follows, the notation

$\bullet$ $U_\Phi:\h\longrightarrow\ell^2$ is used to denote the
{\it analysis operator} of the Bessel sequence $\Phi$ and defined
by $U_\Phi(f):=(\big<f,\varphi_n\big>)_n$ for all $f\in\h$,

$\bullet$ $T_\Phi:\ell^2\longrightarrow\h$ is used to denote the
{\it synthesis operator} of the Bessel sequence $\Phi$ and defined
by $T_\Phi((c_n)_n):=\sum_{n=1}^\infty c_n\varphi_n$ for all
$f\in\h$,

$\bullet$ $S_\Phi:\h\longrightarrow\h$ is used to denote the {\it
frame operator} of the Bessel sequence $\Phi$ and defined by
$S_\Phi(f):=\sum_{n=1}^\infty\big<f,\varphi_n\big>\varphi_n$ for
all $f\in\h$.

Note that $T_\Phi^*=U_\Phi$, where $T_\Phi^*$ is the adjoint of
the operator $T_\Phi$.  Furthermore,  if $\Phi$ is a frame, then
$S_\Phi$ is a bounded, invertible, self-adjoint and positive
operator such that any $f\in\h$ can be expressed as
\begin{equation}\label{02}
f=\sum_{n=1}^\infty\big<f,S_\Phi^{-1}(\varphi_n)\big>\varphi_n
=\sum_{n=1}^\infty\big<f,\varphi_n\big>S_\Phi^{-1}(\varphi_n);
\end{equation}
From now on, the notation $\widetilde{\Phi}$ is used to denote the
sequence $(\widetilde{\varphi}_n)_n$ which defined by
$\widetilde{\varphi}_n:=S_\Phi^{-1}(\varphi_n)$ for all
$n\in{\Bbb N}$, and it is called the {\it canonical dual frame} of
$\Phi$. Moreover, for the frame $\Phi$ which is not a Riesz basis,
there exists infinitely many sequences $\Phi^d=(\varphi_n^d)_n$
such that the following reconstruction formula is hold for all
$f\in\h$
\begin{equation}\label{03}
f=\sum_{n=1}^\infty\big<f,\varphi_n^d\big>\varphi_n
=\sum_{n=1}^\infty\big<f,\varphi_n\big>\varphi_n^d,
\end{equation}
see Theorem 5.2.3 of \cite{c}. Recall that the sequence $\Phi^d$
which satisfies the inequality (\ref{03}) is called the {\it dual
frame} of $\Phi$. In terms of the operators $T_{\Phi}$ and
$U_{\Phi^d}$, the equality (\ref{03}) means that
$T_{\Phi}U_{\Phi^d}=Id_\h=T_{\Phi^d}U_{\Phi}$, where here and in
the sequel, $Id_\h$ is the identity operator on $\h$.

Approximately dual frames as an applicable and interesting duality
principle in the theory of frames was introduced by Christensen
and Laugesen in \cite{app}. Here it should be noted that, the idea
of approximately dual frames can be found in Gabor systems
\cite{hans1,a1,a3}, wavelets \cite{app1,a2,wave1}, coorbit
theory \cite{hans2} and sensor modeling \cite{li}, 
before they introduce it. Moreover, their paper initiated a series of
subsequent publications (see for example
\cite{app1,app2,app3,app4,app5,app7,app8}) and has had a great
impact. In fact, they said that two Bessel sequences
$\Phi=(\varphi_n)_n$ and $\Phi^{ad}=(\varphi_n^{ad})_n$ are
approximately dual frames if $\|Id_\h-T_\Phi U_{\Phi^{ad}}\|<1$ or
$\|Id_\h-T_{\Phi^{ad}} U_\Phi\|<1$. It follows that, the operator
$\A=T_\Phi U_{\Phi^{ad}}$ is invertible and any $f\in\h$ can be
expressed as
\begin{equation}\label{50}
f=\sum_{n=1}^\infty\big<\A^{-1} f,\varphi_n^{ad}\big>\varphi_n.
\end{equation}

Recently, another generalization of duality principle in the
theory of frames has been proposed by Dehghan and Hasankhani-Fard
\cite{app2}. Let us recall from \cite{app2} that, the frame
$\Phi^{gd}=(\varphi_n^{gd})_n$ is a generalized dual frame or
g-dual frame of $\Phi$ with corresponding invertible operator (or
with invertible operator) $\A\in B(\h)$, if we have the following
inequality for all $f\in\h$
\begin{equation}\label{04}
f=\sum_{n=1}^\infty\big<\A f,\varphi_n^{gd}\big>\varphi_n.
\end{equation}
In terms of the operators $T_{\Phi}$ and $U_{\Phi^{gd}}$, the
equality (\ref{04}) means that $T_{\Phi}U_{\Phi^{gd}}=\A^{-1}$. It
follows that the operator $\mathcal A$ in (\ref{04}) is unique.
Moreover, a simple observation shows that if two frames are
approximately dual frames, then they are g-dual frames. Although
this is not unexpected, Example 4.1 of \cite{app2} illustrated
that the set of approximately duals of a frame is a proper subset
of the set of its g-duals. Moreover, recall from \cite[Theorem
3.1]{app2} that the set of all g-dual frames of the frame $\Phi$
with corresponding invertible operator $\A^* S_\Phi^{-1}$ are
precisely the sequences of the from
$$\Phi^{gd}={\bigg(}\A \varphi_n+\psi_n-\sum_{m=1}^\infty
\big<S_\Phi^{-1}\varphi_n,\varphi_m\big>\psi_m{\bigg)}_n,$$ where
$\Psi=(\psi_n)_n$ is a Bessel sequence in $\h$ and $\A\in B(\h)$
is an invertible operator; In particular $T_\Phi U_\Psi=\A^*
S_\Phi^{-1}$.

In the present work, we introduce a general method of constructing
approximately duals (resp. g-duals) frame for a given frame and
given operator by using of bounded operators which our explicit
construction can be easily applied for Gabor frames.  We then show
that for a perturbed frame, one can construct always an
approximately dual (resp. g-dual) frame which is close to the
approximately dual (resp. g-dual) frame of the original frame.
Finally, we show that by choosing an appropriate dual generator
for a Gabor frame provides more precise results.


\section{\large\bf The basic results}

We commence this section with the following result which gives a
sufficient and necessary condition for two frames $\Phi$ and
$\Psi$ under which they are g-dual frames. To this end, recall
that every bounded and positive operator ${\cal
C}:\h\rightarrow\h$ has a unique bounded and positive square root
${\cal C}^{\frac{1}{2}}$. Moreover, if the operator $\cal C$ is
self-adjoint (resp. invertible), then ${\cal C}^{\frac{1}{2}}$ is
also self-adjoint (resp. invertible), see for example Lemma 2.4.4
of \cite{c}.

\begin{lemma}\label{appdual}
Let $\Phi$ and $\Psi$ be two frames for $\h$ with Bessel bounds
$\M_\Phi$ and $\M_\Psi$, respectively. Then $\Psi$ is a g-dual
frame of $\Phi$ if and only if there exists an invertible operator
${\cal D}\in B(\h)$ such that $T_\Phi
U_{\Psi}=S_\Phi^{\frac{1}{2}}{\cal D}$ and ${\cal D}{\cal
D}^*\leq\M_\Psi\;Id_\h$.
\end{lemma}
{\noindent Proof.} The proof of the ``only if" part is trivial. To
prove the ``if" part, suppose that $\Psi$ is a g-dual frame of
$\Phi$. Observe that
\begin{equation}\label{06}
\frac{1}{\M_{\Psi}}\;\|T_{\Psi} U_\Phi
f\|^2\leq\sum_{n=1}^\infty|\big<f,\varphi_n\big>|^2,
\end{equation}
for all $f\in\h$; Indeed,
\begin{eqnarray*}
\|T_{\Psi} U_\Phi f\|&=&\sup_{\|g\|=1}|\big<T_{\Psi} U_\Phi f,g\big>|\\
&=&\sup_{\|g\|=1}|\big<f,T_\Phi U_{\Psi} g\big>|\\
&=&\sup_{\|g\|=1}{\bigg|}\big<f,\sum_{n=1}^\infty\big<g,\psi_n
\big>\varphi_n\big>{\bigg|}\\
&=&\sup_{\|g\|=1}{\bigg|}\sum_{n=1}^\infty\overline{\big<g,\psi_n
\big>}\;\big<f,\varphi_n\big>{\bigg|}\\
&\leq&\sup_{\|g\|=1}{\bigg(}\sum_{n=1}^\infty|\big<g,\psi_n
\big>|^2{\bigg)}^{\frac{1}{2}}
{\bigg(}\sum_{n=1}^\infty|\big<f,\varphi_n\big>|^2{\bigg)}^{\frac{1}{2}}\\
&\leq&\sup_{\|g\|=1}\sqrt{\M_{\Psi}}\;\|g\|{\bigg(}\sum_{n=1}^\infty|
\big<f,\varphi_n\big>|^2{\bigg)}^{\frac{1}{2}}\\
&\leq&\sqrt{\M_{\Psi}}\;{\bigg(}\sum_{n=1}^\infty|
\big<f,\varphi_n\big>|^2{\bigg)}^{\frac{1}{2}}.
\end{eqnarray*}
Note that, in terms of the inner product of $\h$, the equality
(\ref{06}) means that
$${\bigg<}T_{\Psi} U_\Phi f,T_{\Psi} U_\Phi
f{\bigg>}\leq M_{\Psi}\big<S_\Phi f,f\big>,\quad\quad\quad
(f\in\h).$$ Therefore, we have the following inequality
\begin{equation}\label{20}
{\bigg<}(T_\Phi U_{\Psi})(T_\Phi U_{\Psi})^* f,f{\bigg>}\leq
M_{\Psi}{\bigg<}S_\Phi^{\frac{1}{2}}S_\Phi^{\frac{1}{2}}
f,f{\bigg>},
\end{equation}
for all $f\in\h$. We now invoke Theorem 1 of \cite{do} and the
invertibility of the operator $T_\Phi U_{\Psi}$, to conclude that
there exists an invertible operator ${\mathcal D}\in B({\mathcal
H})$ such that $T_\Phi U_{\Psi}=S_\Phi^{\frac{1}{2}}{\cal D}$.
Moreover, inequality (\ref{20}) means that
$$S_\Phi^{\frac{1}{2}}{\cal D}{\cal D}^*S_\Phi^{\frac{1}{2}}=
(S_\Phi^{\frac{1}{2}}{\cal D})(S_\Phi^{\frac{1}{2}}{\cal D})^*\leq
\M_{\Psi}S_\Phi^{\frac{1}{2}}S_\Phi^{\frac{1}{2}}.$$ Therefore,
\cite[Theorem 2.2.5(2)]{mur} implies that ${\cal D}{\cal D}^*\leq
\M_{\Psi}\;Id_\h$.$\e$\\


By the same argument as above, we can prove the following result.

\begin{lemma}\label{appdual2}
Let $\Phi$ and $\Psi$ be two frames for $\h$ with Bessel bounds
$\M_\Phi$ and $\M_\Psi$, respectively. Then $\Psi$ is an
approximately dual frame of $\Phi$ if and only if there exists an
operator ${\mathcal D}\in B({\mathcal H})$ such that $T_\Phi
U_{\Psi}=S_\Phi^{\frac{1}{2}}{\cal D}$, ${\cal D}{\cal
D}^*\leq\M_\Psi\;Id_\h$ and
$\|Id_\h-S_{\Phi}^{\frac{1}{2}}{\mathcal D}\|<1$.\\
\end{lemma}


From now on, for the frame $\Phi=(\varphi_n)_n$, we assume that
$\verb"ran"_{B(\h,\ell^2)}(T_\Phi)$ stands for the set of all
right annihilator of the operator $T_\Phi$ in $B(\h,\ell^2)$; That
is,
$$\verb"ran"_{B(\h,\ell^2)}(T_\Phi):={\bigg\{}\Theta\in B(\h,\ell^2):~
T_\Phi\Theta=0{\bigg\}}.$$ By the use of \cite[Theorem 5.2.2]{c}
and a routine computation, we can see that the frame $\Phi$ is a
Riesz basis for $\h$ if and only if
$\verb"ran"_{B(\h,\ell^2)}(T_\Phi)=\{0\}$.

The following theorem gives a new characterization of all
approximately dual frames of $\Phi=(\varphi_n)_n$ in terms of the
operators in $\verb"ran"_{B(\h,\ell^2)}(T_\Phi)$; In fact, this is
a construction result for approximately dual frames associated
with a given frame. For the rest of the paper the letter
``$\big<\cdot,\cdot\big>_{\ell^2}$" means the inner product of
$\ell^2$.

\begin{theorem}\label{70}
Let $\Phi=(\varphi_n)_n$ be a frame for $\h$. Then, the
approximately dual frames of $\Phi$ are precisely the families
$\Phi^{ad}=(\varphi_n^{ad})_n$ such that
$$\varphi^{ad}_n={\mathcal D}^*S_{\Phi}^{-{\frac{1}{2}}}\varphi_n+
\Theta^*(\delta_n),\quad\quad(n\in{\Bbb N})$$ where $\Theta\in
\verb"ran"_{B(\h,\ell^2)}(T_\Phi)$ and $\mathcal D$ is a bounded
operator on $\h$ for which
$\|Id_\h-S_{\Phi}^{\frac{1}{2}}{\mathcal D}\|<1$. In particular,
$T_\Phi U_{\Phi^{ad}}=S_\Phi^{\frac{1}{2}}{\cal D}$.
\end{theorem}
{\noindent Proof.} First, assume that
$\Phi^{ad}=(\varphi^{ad}_n)_n$ is a sequence in $\h$ such that
$$\varphi^{ad}_n={\mathcal D}^*S_{\Phi}^{-{\frac{1}{2}}}\varphi_n+
\Theta^*(\delta_n),\quad\quad\quad(n\in{\Bbb N})$$ where
$\Theta\in \verb"ran"_{B(\h,\ell^2)}(T_\Phi)$ and $\mathcal D$ is
a bounded operator on $\h$ for which
$\|Id_\h-S_{\Phi}^{\frac{1}{2}}{\mathcal D}\|<1$. Now, since the
sequence
$\Phi^{-{\frac{1}{2}}}:=(S_{\Phi}^{-{\frac{1}{2}}}\varphi_n)_n$ is
a tight frame with frame bound equal to $1$, we can see easily
that $\Phi^{ad}$ is a Bessel sequence with Bessel bound
$\M_{\Phi^{ad}}=\|{\mathcal D}\|^2+\|\Theta\|^2+2\|{\mathcal
D}\|\|\Theta\|$. We also observed that
\begin{eqnarray*}
T_\Phi U_{\Phi^{ad}} f=S_\Phi^{{\frac{1}{2}}}{\mathcal D};
\end{eqnarray*}
This is because of,
\begin{eqnarray*}
U_{\Phi^{ad}}(f)={\bigg(}\big<f,\varphi^{ad}_n\big>{\bigg)}_n&=&{\bigg(}\big<f,{\mathcal
D}^*S_{\Phi}^{-{\frac{1}{2}}}\varphi_n\big>+\big<f,\Theta^*(\delta_n)\big>{\bigg)}_n\\
&=&{\bigg(}\big<S_{\Phi}^{-{\frac{1}{2}}}{\mathcal
D}f,\varphi_n\big>+\big<\Theta f,\delta_n\big>_{\ell^2}{\bigg)}_n\\
&=&{\bigg(}\big<S_{\Phi}^{-{\frac{1}{2}}}{\mathcal
D}f,\varphi_n\big>+(\Theta f)_n{\bigg)}_n\\
&=&U_\Phi S_\Phi^{-{\frac{1}{2}}}{\mathcal D}f+\Theta f,
\end{eqnarray*}
for all $f\in\h$, where $(\Theta f)_n$ is the nth term of the
sequence $\Theta f$. Moreover,
\begin{eqnarray*}
\big<{\mathcal D}{\mathcal D}^*f,f\big>=\|{\mathcal
D}^*f\|^2\leq\|{\mathcal
D}\|^2\big<f,f\big>\leq\M_{\Phi^{ad}}\big<f,f\big>,
\end{eqnarray*}
for all $f\in\h$. We now invoke Lemma \ref{appdual2} to conclude
that $\Phi^{ad}$ is an approximately dual frame of $\Phi$.

Conversely, let $\Phi^{ad}=(\varphi_n^{ad})_n$ be an approximately
dual frame of $\Phi=(\varphi_n)_n$ with Bessel bound
$\M_{\Phi^{ad}}$. Then, in view of Lemma \ref{appdual2}, there
exists an invertible operator $\mathcal D$ in $B(\mathcal H)$ such
that $T_\Phi U_{\Phi^{ad}}=S_\Phi^{\frac{1}{2}}{\mathcal D}$,
${\cal D}{\cal D}^*\leq\M_{\Phi^{ad}}\;Id_\h$ and
$$\|Id_\h-S_{\Phi}^{\frac{1}{2}}{\mathcal D}\|=\|Id_\h-T_\Phi U_{\Phi^{ad}}\|<1.$$ Now,
define the operator $\Theta$ from $\h$ into $\ell^2$ by
$$\Theta:=U_{\Phi^{ad}}-U_\Phi S_\Phi^{-{\frac{1}{2}}} {\mathcal D}.$$
Observe that the operator $\Theta$ is in
$\verb"ran"_{B(\h,\ell^2)}(T_\Phi)$. Moreover, for each $n\in{\Bbb
N}$ we have
\begin{eqnarray*}
\Theta^*(\delta_n)&=&T_{\Phi^{ad}}(\delta_n)-{\mathcal
D}^*S_\Phi^{-{\frac{1}{2}}}T_\Phi(\delta_n)\\
&=&\varphi^{ad}_n-{\mathcal D}^*S_\Phi^{-{\frac{1}{2}}}\varphi_n,
\end{eqnarray*}
and this completes the proof.$\e$\\


As applications of Theorem \ref{70}, we have the following two
results which gives a construction result for approximately dual
frames associated with a given frame and given operator.

\begin{theorem}\label{4000}
Let $\Phi=(\varphi_n)_n$ be a frame for $\h$ with Bessel bound
$\M_\Phi$, and let $\mathcal D$ be an operator in $B(\h)$ such
that $\|S_{\Phi}^{-{\frac{1}{2}}}-{\mathcal
D}\|<\frac{1}{\sqrt{\M_\Phi}}$. Then, the sequence
$$\Phi^{ad}={\bigg(}{\mathcal D}^*S_{\Phi}^{-{\frac{1}{2}}}\varphi_n+
\Theta^*(\delta_n){\bigg)}_n$$ is an approximately dual frame of
$\Phi$ for which $T_\Phi U_{\Phi^{ad}}=S_\Phi^{\frac{1}{2}}{\cal
D}$, where $\Theta$ is an arbitrary element of
$\verb"ran"_{B(\h,\ell^2)}(T_\Phi)$.
\end{theorem}
{\noindent Proof.} Note that
\begin{eqnarray*}
\|Id_\h-S_{\Phi}^{\frac{1}{2}}{\mathcal D}\|&=&
\|S_{\Phi}^{\frac{1}{2}}S_{\Phi}^{-\frac{1}{2}}-
S_{\Phi}^{\frac{1}{2}}{\mathcal D}\|\\
&\leq&\|S_{\Phi}^{\frac{1}{2}}\|\;\|S_{\Phi}^{-\frac{1}{2}}-{\mathcal
D}\|\\
&=&\sqrt{\|(S_\Phi^{\frac{1}{2}})^*S_\Phi^{\frac{1}{2}}\|}
\;\|S_{\Phi}^{-\frac{1}{2}}-{\mathcal
D}\|\\
&=&\sqrt{\|S_\Phi\|}\;\|S_{\Phi}^{-\frac{1}{2}}-{\mathcal D}\|\\
&=&\sqrt{\M_\Phi}\;\|S_{\Phi}^{-\frac{1}{2}}-{\mathcal D}\|.
\end{eqnarray*}
So, the proof will follow from Theorem \ref{70}.$\e$\\


\begin{theorem}\label{40}
Let $\Phi=(\varphi_n)_n$ be a frame for $\h$, and let $\mathcal A$
be an operator in $B(\h)$ such that $\|Id_\h-{\mathcal A}\|<1$.
Then, the Bessel sequence $\Phi^{ad}=(\varphi_n^{ad})_n$ is an
approximately dual frame of $\Phi$ for which $T_\Phi
U_{\Phi^{ad}}={\mathcal A}$ if and only if there exists $\Theta\in
\verb"ran"_{B(\h,\ell^2)}(T_\Phi)$ such that
$\varphi^{ad}_n={\mathcal A}^*\widetilde{\varphi}_n+
\Theta^*(\delta_n)$ for all $n\in{\Bbb N}$.
\end{theorem}


The following two corollaries illustrate the construction given in
Theorems \ref{4000} and \ref{40}, which are very applicable for
the construction of approximately dual Gabor frames with a desired
approximation rate. Here it should be noted that, if $\Phi$ and
$\Psi$ are approximately dual frames, then we say that they have
approximation rate $0<\varepsilon<1$, if $\|Id_\h-T_\Phi
U_\Psi\|\leq\varepsilon$.

\begin{corollary}\label{4001}
Let $\Phi=(\varphi_n)_n$ be a frame for $\h$ with Bessel bound
$\M_\Phi$ and a dual frame $\Phi^d=(\varphi_n^d)_n$. Then, for
each operator $\mathcal D$ in $B(\h)$ such that
$\|S_{\Phi}^{-{\frac{1}{2}}}-{\mathcal
D}\|<\frac{1}{\sqrt{\M_\Phi}}$, the sequence
$$\Phi^{ad}={\bigg(}{\mathcal D}^*S_{\Phi}^{-{\frac{1}{2}}}\varphi_n
-\varphi_n+S_\Phi(\varphi_n^d){\bigg)}_n$$ is an approximately
dual frame of $\Phi$ for which $T_\Phi
U_{\Phi^{ad}}=S_\Phi^{\frac{1}{2}}{\cal D}$.\\
\end{corollary}

\begin{corollary}\label{200}
Let $\Phi=(\varphi_n)_n$ be a frame for $\h$ with a dual frame
$\Phi^d=(\varphi_n^d)_n$. Let also $\A$ be a bounded operator on
$\h$ such that $\|Id_\h-{\mathcal A}\|<1$. Then, the sequence
$\Phi^{ad}=(\varphi_n^{ad})_n$ defined by
$$\varphi^{ad}_n={\mathcal
A}^*\widetilde{\varphi}_n-\varphi_n+S_\Phi(\varphi_n^d),
\quad\quad\quad(n\in{\Bbb N})$$ is an approximately dual frame of
$\Phi$ for which $T_\Phi U_{{\Phi}^{ad}}={\mathcal A}$.
\end{corollary}


In the following, we use the perturbation idea to construct
approximately dual frames. The notation $\Phi-\Psi$ (resp.
$\Phi^{ad}-\Psi^{ad}$) in this theorem denotes the sequence
$(\varphi_n-\psi_n)_n$ (resp. $(\varphi_n^{ad}-\psi_n^{ad})_n$).

\begin{theorem}\label{80}
Let $\Phi=(\varphi_n)_n$ and $\Psi=(\psi_n)_n$ be two frames for
$\h$ and let $\Phi^{ad}=(\varphi_n^{ad})_n$ be a fixed
approximately dual frame of $\Phi$. If $\Phi-\Psi$ is a Bessel
sequence with sufficiently small bound $\M_{\Phi-\Psi}$, then
there exists an approximately dual frame $\Psi^{ad}=(\psi_n^{ad})$
of $\Psi$ such that $\Phi^{ad}-\Psi^{ad}$ is also a Bessel
sequence and its bound is a multiple of $\M_{\Phi-\Psi}$. In
particular, $T_\Phi U_{\Phi^{ad}}=T_\Psi U_{\Psi^{ad}}$.
\end{theorem}
{\noindent Proof}. In view of Theorems \ref{70} and \ref{40},
there exists an operator $\mathcal A$ in $B(\mathcal H)$ and
$\Theta\in \verb"ran"_{B(\h,\ell^2)}(T_\Phi)$ such that

{(a)} $\|Id_\h-{\mathcal A}\|<1$,

(b) $\varphi^{ad}_n={\mathcal A}^*\widetilde{\varphi}_n+
\Theta^*(\delta_n),$ for all $n\in{\Bbb N}$,

(c) $T_\Phi U_{\Phi^{ad}}=\A$.\\
Now, let $\Omega:=(\omega_n)_n$ be a sequence in $\h$ such that
for each $n\in{\Bbb N}$
$$\omega_n={\mathcal A}^*\widetilde{\psi}_n+
\Theta^*(\delta_n).$$ A routine computation shows that the
sequence $\Omega$ is a Bessel sequence. Moreover, for each
$f\in\h$ we have
\begin{eqnarray*}
{\mathcal A}^{-1}T_\Psi U_\Omega f&=&{\mathcal
A}^{-1}{\bigg(}\sum_{n=1}^\infty\big<f,\omega_n\big>\psi_n{\bigg)}\\
&=&{\mathcal A}^{-1}{\bigg(}\sum_{n=1}^\infty\big<f,{\mathcal
A}^*\widetilde{\psi}_n\big>\psi_n
+\sum_{n=1}^\infty\big<f,\Theta^*(\delta_n)\big>\psi_n{\bigg)}\\
&=&{\mathcal A}^{-1}{\bigg(}\sum_{n=1}^\infty\big<{\mathcal
A}f,\widetilde{\psi}_n\big>\psi_n
+\sum_{n=1}^\infty\big<\Theta f,\delta_n\big>_{\ell^2}\psi_n{\bigg)}\\
&=&{\mathcal A}^{-1}{\bigg(}{\mathcal
A}f+T_\Psi\Theta f{\bigg)}\\
&=&f+{\mathcal A}^{-1}T_\Psi\Theta f.
\end{eqnarray*}
Therefore,
\begin{eqnarray*}
\|f-{\mathcal A}^{-1}T_\Psi U_\Omega f\|= \|{\mathcal
A}^{-1}T_\Psi\Theta f\|&=&\|{\mathcal A}^{-1}T_\Psi\Theta
f-{\mathcal A}^{-1}T_\Phi\Theta f\|\\
&\leq &\sqrt{\M_{\Phi-\Psi}}\;\|\Theta\|\;\|\A^{-1}\|\;\|f\|,
\end{eqnarray*}
for all $f\in\h$. Hence, the operator ${\mathcal C}:=T_\Omega
U_\Psi {{\mathcal A}^*}^{-1}$ is invertible for sufficiently small
$\M_{\Phi-\Psi}>0$. In particular, the operator $T_{\Psi}U_\Omega$
is invertible and for each $f\in\h$ we have
\begin{eqnarray*}
f=(T_{\Psi}U_\Omega)(T_{\Psi}U_\Omega)^{-1}f&=&
\sum_{n=1}^\infty\big<(T_{\Psi}U_\Omega)^{-1}f,\omega_n\big>\psi_n\\
&=&\sum_{n=1}^\infty\big<f,(T_\Omega U_\Psi)^{-1}\omega_n\big>\psi_n\\
&=&\sum_{n=1}^\infty\big<{\mathcal A}^{-1}f,(T_\Omega
U_\Psi{{\mathcal
A}^*}^{-1})^{-1}\omega_n\big>\psi_n\\
&=&\sum_{n=1}^\infty\big<{\mathcal A}^{-1}f,{\mathcal
C}^{-1}\omega_n\big>\psi_n.
\end{eqnarray*}
It follows that $\Psi^{ad}:=({\cal C}^{-1} \omega_n)_n$ is a frame
for $\h$ such that $T_\Psi U_{\Psi^{ad}}={\mathcal A}$ and
therefore $\Psi^{ad}$ is an approximately dual frame of $\Psi$.

Now, the proof will be completed if we show that
$\Phi^{ad}-\Psi^{ad}$ is also a Bessel sequence and its bound is a
multiple of $\M_{\Phi-\Psi}$. To this end, we note that
\begin{equation}\label{07}
\|{\cal C}^{-1}\|\leq\frac{1}{1-\|Id_\h-{\cal C}\|} \leq\frac{1}
{1-\sqrt{\M_{\Phi-\Psi}}\;\|\Theta\|\|\A^{-1}\|},
\end{equation}
\begin{equation}\label{08}
\|Id_\h-{\cal C}^{-1}\|\leq\|{\cal C}^{-1}\|\|Id_\h-{\cal C}\|\leq
\frac{\sqrt{\M_{\Phi-\Psi}}\;\|\Theta\|\|\A^{-1}\|}
{1-\sqrt{\M_{\Phi-\Psi}}\;\|\Theta\|\|\A^{-1}\|},
\end{equation}
and
\begin{eqnarray*}\label{09}
\|S_\Phi-S_\Psi\|&=&\|T_\Phi U_\Phi-T_\Phi U_\Psi+T_\Phi U_\Psi-T_\Psi U_\Psi\|\\
&\leq&\|T_\Phi-T_\Psi\|(\|T_\Phi\|+\|T_\Psi\|)\\
&\leq&\sqrt{\M_{\Phi-\Psi}}\;(\sqrt{M_\Phi}+\sqrt{M_\Psi}\;)
\end{eqnarray*}
Moreover, for
all $(c_n)_n\in\ell^2$ and all $f\in{\mathcal H}$ we have
\begin{eqnarray*}
\big<T_{\Phi^{ad}}((c_n)_n),f\big>-\big<T_\Omega((c_n)_n),f\big>&=&
\big<(c_n)_n,U_{\Phi^{ad}} f\big>_{\ell^2}-\big<(c_n)_n,U_\Omega f\big>_{\ell^2}\\
&=&\big<(c_n)_n,U_{\Phi}S_\Phi^{-1}\A f\big>_{\ell^2}-
\big<(c_n)_n,U_{\Psi}S_\Psi^{-1}\A f\big>_{\ell^2}\\
&=&\big<T_{\Phi}((c_n)_n),S_\Phi^{-1}\A f\big>-\big<T_{\Psi}((c_n)_n),S_\Psi^{-1}\A f\big>\\
&=&\big<T_{\Phi}((c_n)_n),(S_\Phi^{-1}-S_\Psi^{-1})\A f\big>\\&&+
\big<(T_{\Phi}-T_{\Psi})((c_n)_n),S_\Psi^{-1}\A f\big>.
\end{eqnarray*}
Hence,
\begin{eqnarray*}
\|T_{\Phi^{ad}}((c_n)_n)-T_\Omega((c_n)_n)\|&=&\sup_{\|g\|=1}{\bigg|}\big<
T_{\Phi^{ad}}((c_n)_n)-T_\Omega((c_n)_n),g\big>{\bigg|}\\
&=&\sup_{\|g\|=1}{\bigg|}\big<T_{\Phi}((c_n)_n),(S_\Phi^{-1}-S_\Psi^{-1})\A
g\big>\\&&+
\big<T_{\Phi}((c_n)_n)-T_{\Psi}((c_n)_n),S_\Psi^{-1}\A g\big>{\bigg|}\\
&\leq&\|S_\Phi^{-1}-S_\Psi^{-1}\|\;\|\A\|\;\|T_{\Phi}((c_n)_n)\|\\&&+
\|S_\Psi^{-1}\|\;\|\A\|\;\|T_{\Phi}((c_n)_n)-T_{\Psi}((c_n)_n)\|\\
&\leq&\|S_\Psi^{-1}\|\;\|S_\Phi-S_\Psi\|\;\|S_\Phi^{-1}\|\;\|\A\|\;\|T_{\Phi}((c_n)_n)\|\\&&+
\|S_\Psi^{-1}\|\;\|\A\|\;\|T_{\Phi}((c_n)_n)-T_{\Psi}((c_n)_n)\|\\
&\leq&\sqrt{\M_{\Phi-\Psi}}\;{\bigg(}\frac{1}{m_\Phi
m_\Psi}(M_\Phi+\sqrt{M_\Psi
M_\Phi}\;)\|\A\|+\frac{1}{m_\Psi}\|\A\| {\bigg)},
\end{eqnarray*}
where, $\m_\Phi, \M_\Phi$ and $\m_\Psi, \M_\Psi$ are the bounds of
the frames $\Phi$ and $\Psi$, respectively. It follows that the
sequence $\Phi^{ad}-\Psi^{ad}$ is a Bessel sequence with Bessel
bound
$$\M_{\Phi^{ad}-\Psi^{ad}}:=\M_{\Phi-\Psi}\;{\bigg(}\frac{\|\A^{-1}\|}
{1-\sqrt{\M_{\Phi-\Psi}}\;\|\Theta\|\|\A^{-1}\|}{\bigg)}^2
{\bigg(}\|\Theta\|\sqrt{M_{\Phi^{ad}}}+\frac{\|\A\|(m_\Phi+M_\Phi+\sqrt{M_\Psi
M_\Phi}\;)}{m_\Phi m_\Psi}{\bigg)}^2$$ and this completes the
proof; Indeed,
\begin{eqnarray*}
\|T_{\Phi^{ad}}((c_n)_n)-T_{\Psi^{ad}}((c_n)_n)\|&=&
\|T_{\Phi^{ad}}((c_n)_n)-T_{\Psi^{ad}}((c_n)_n)
-{\mathcal C}^{-1}T_{\Phi^{ad}}((c_n)_n)+{\mathcal C}^{-1}T_{\Phi^{ad}}((c_n)_n)\|\\
&=&\|T_{\Phi^{ad}}((c_n)_n)-{\mathcal C}^{-1}T_{\Omega}((c_n)_n)
-{\mathcal C}^{-1}T_{\Phi^{ad}}((c_n)_n)+{\mathcal C}^{-1}T_{\Phi^{ad}}((c_n)_n)\|\\
&\leq&\|Id_{\mathcal H}-{\mathcal
C}^{-1}\|\;\|T_{\Phi^{ad}}((c_n)_n)\| +\|{\mathcal C}^{-1}\|\;
\|T_{\Phi^{ad}}((c_n)_n)-T_{\mathcal F}((c_n)_n)\|
\end{eqnarray*}
for all $(c_n)_n\in\ell^2$.$\e$\\


An argument similar to the proof of Theorems \ref{70}, \ref{40}
and \ref{80} with the aid of Lemma \ref{appdual} gives the
following generalization of that theorems. The details are
omitted.

\begin{theorem}
Let $\Phi=(\varphi_n)_n$ be a frame for $\h$ and let $\mathcal A$
be an invertible operator in $B(\h)$. Then
$\Phi^{gd}=(\varphi_n^{gd})_n$ is a g-dual frame of $\Phi$ with
corresponding invertible operator ${\mathcal A}\in B(\h)$ if and
only if there exists $\Theta\in \verb"ran"_{B(\h,\ell^2)}(T_\Phi)$
such that $\varphi^{gd}_n={{\mathcal
A}^{-1}}^*\widetilde{\varphi}_n+ \Theta^*(\delta_n)$ for all
$n\in{\Bbb N}$.
\end{theorem}


\begin{theorem}\label{pergdd}
Let $\Phi=(\varphi_n)_n$ and $\Psi=(\psi_n)_n$ be two frames for
$\h$ and let $\Phi^{gd}=(\varphi^{gd}_n)_n$ be a fixed g-dual
frame of $\Phi$ with corresponding invertible operator $\mathcal
A$. If $\Phi-\Psi$ is a Bessel sequence with sufficiently small
bound $\M_{\Phi-\Psi}$, then there exists a g-dual frame
$\Psi^{gd}=(\psi^{gd}_n)_n$ of $\Psi$ such that
$\Phi^{gd}-\Psi^{gd}$ is also a Bessel sequence and its bound is a
multiple of $\M_{\Phi-\Psi}$. In particular, $T_\Phi
U_{\Phi^{gd}}={\mathcal A}^{-1}=T_\Psi U_{\Psi^{gd}}$.
\end{theorem}


Since each frame $\Phi$ is a g-dual frame for itself with
corresponding invertible operator $S_\Phi^{-1}$, hence an argument
similar to the proof of Theorem \ref{80} with the aid of Theorem
\ref{pergdd} gives the following result.

\begin{theorem}\label{100}
Let $\Phi=(\varphi_n)_n$ and $\Psi=(\psi_n)_n$ be two frames for
$\h$. Let also $\M_{\Phi-\Psi}$ is the Bessel bound of the Bessel
sequence $\Phi-\Psi$. Then there exists a g-dual frame
$\Psi^{gd}=(\psi^{gd}_n)_n$ of $\Psi$ such that $\Phi-\Psi^{gd}$
is also a Bessel sequence, its bound is a multiple of
$\M_{\Phi-\Psi}$ and $T_\Psi U_{\Psi^{gd}}=S_\Phi$.
\end{theorem}


Now, let $\Phi=(\varphi_n)_n$ and $\Psi=(\psi_n)_n$ are Riesz
bases with Bessel bounds $\M_\Phi$ and $\M_\Psi$, respectively.
Then, there exists an invertible operator $\cal D$ on $\h$ such
that ${\cal D}\varphi_n=\psi_n$ for all $n\in{\Bbb N}$. Observe
that the sequence $\Phi-\Psi$ is a Bessel sequence with Bessel
bound
$$\M_{\Phi-\Psi}=\min{\bigg\{}\M_\Phi\|Id_\h-{\cal D}\|^2,
\M_\Psi\|Id_\h-{\cal D}^{-1}\|^2{\bigg\}};$$ Indeed, if $(c_n)_n$
is an arbitrary element of $\ell^2$, then on the one hand we have
\begin{eqnarray*}
\|(T_\Phi-T_\Psi)((c_n)_n)\|&=&{\bigg\|}\sum_{n=1}^\infty
c_n(\varphi_n-{\cal D}(\varphi_n)){\bigg\|}\\
&\leq&\|Id_\h-{\cal D}\|\;\|T_\Phi((c_n)_n)\|\\
&\leq&\sqrt{\M_\Phi}\;\|Id_\h-{\cal D}\|\;\|(c_n)_n\|_2,
\end{eqnarray*}
and on the other hand, by the same argument we can see that
\begin{eqnarray*}
\|(T_\Phi-T_\Psi)((c_n)_n)\|\leq\sqrt{\M_\Psi}\;\|Id_\h-{\cal
D}^{-1}\|\;\|(c_n)_n\|_2.
\end{eqnarray*}
Hence as an immediate corollary from Theorem \ref{100}, we have
the following result.

\begin{theorem}\label{110}
Let $\Phi=(\varphi_n)_n$ and $\Psi=(\psi_n)_n$ be two Riesz bases
for $\h$ with Bessel bounds $\M_\Phi$ and $\M_\Psi$, respectively.
Then $\Phi-\Psi$ is a Bessel sequence with Bessel bound
$\M_{\Phi-\Psi}=\min\{\M_\Phi\|Id_\h-{\cal D}\|^2,
\M_\Psi\|Id_\h-{\cal D}^{-1}\|^2\}$, where $\cal D$ is an operator
such that ${\cal D}\varphi_n=\psi_n$ for all $n\in{\Bbb N}$.
Moreover, there exists a g-dual frame $\Psi^{gd}=(\psi_n^{gd})_n$
of $\Psi$ such that $\Phi-\Psi^{gd}$ is also a Bessel sequence,
its bound is a multiple of $\M_{\Phi-\Psi}$ and $T_\Psi
U_{\Psi^{gd}}=S_\Phi$. In particular, $\Psi^{gd}$ is a Riesz
basis.
\end{theorem}


Following Balan \cite{balan}, we say that two frames $\Phi$ and
$\Psi$ for $\h$ are equivalent frames, if there exists a bounded
invertible operator $Q:{\h}\rightarrow{\h}$ such that
$\varphi_n=Q(\psi_n)$ for all $n\in{\Bbb N}$ or equivalently ${\rm
Range}(U_\Phi)={\rm Range}(U_\Psi)$. Furthermore, the frame $\Psi$
is partial equivalent with the frame $\Phi$, if there exists a
bounded operator $Q:{\h}\rightarrow{\h}$ (not necessarily
invertible) such that $\varphi_n=Q(\psi_n)$ for all $n\in{\Bbb N}$
or equivalently ${\rm Range}(U_\Phi)\subsetneq {\rm
Range}(U_\Psi)$. Observe that, in the case where ${\rm
Range}(U_\Phi)\subseteq {\rm Range}(U_\Psi)$, then
$T_{\widetilde{\Psi}}U_{\widetilde{\Phi}}$ is a right inverse of
$T_\Phi U_\Psi$; Indeed, if $f$ is an arbitrary element of $\h$.
Then, there exists $h\in\h$ such that $U_\Phi S_\Phi^{-1} f=U_\Psi
h$, and therefore
\begin{eqnarray*}
T_\Phi U_\Psi T_{\widetilde{\Psi}}U_{\widetilde{\Phi}}(f)&=&
T_\Phi U_\Psi S_\Psi^{-1}T_{{\Psi}}U_{{\Phi}}S_\Phi^{-1}(f)\\
&=&T_\Phi U_\Psi S_\Psi^{-1}T_{{\Psi}}U_{{\Psi}}(h)\\
&=&T_\Phi U_\Psi(h)\\
&=&T_\Phi U_\Phi S_\Phi^{-1} f=f.
\end{eqnarray*}
So we have the following result for equivalent frames.

\begin{proposition}
Let $\Phi$ and $\Psi$ be two frames for $\cal H$ such that ${\rm
Range}(U_\Phi)={\rm Range}(U_\Psi)$. Then $\Phi$ and $\Psi$ are
g-dual frames. In particular, $(T_\Phi U_\Psi)^{-1}=
T_{\widetilde{\Psi}}U_{\widetilde{\Phi}}$.
\end{proposition}


We conclude this section by the following result which is of
interest in its own right.

\begin{proposition}\label{equiv}
Let $\Phi$ and $\Psi$ be two frames for $\cal H$ such that ${\rm
Range}(U_\Phi)\subsetneq {\rm Range}(U_\Psi)$ or ${\rm
Range}(U_\Psi)\subsetneq {\rm Range}(U_\Phi)$. Then $\Phi$ and
$\Psi$ are not g-dual frames. In particular, they are not
approximately dual frames.
\end{proposition}
{\noindent Proof.} Let $\Phi$ and $\Psi$ be two frames for $\cal
H$ such that ${\rm Range}(U_\Phi)\subsetneq {\rm Range}(U_\Psi)$.
As we saw in above discussion, this condition implies that
$T_{\widetilde{\Psi}}U_{\widetilde{\Phi}}$ is a right inverse of
$T_\Phi U_\Psi$. Now we show that it is not a left inverse for
$T_\Phi U_\Psi$. To this end, assume that $(c_n)_n$ is an element
of ${\rm Range}(U_\Psi)$ such that $(c_n)_n\notin {\rm
Range}(U_\Phi)$. Then, there exists $g_1\in\h$ and $0\neq
(d_n)_n\in\ker(T_\Phi)$ such that $(c_n)_n=U_\Phi g_1+(d_n)_n$;
This is because of $\ell^2={\rm Range}(U_\Phi)\oplus\ker(T_\Phi)$.
It follows that $(d_n)_n=U_\Phi g_1-(c_n)_n\in {\rm
Range}(U_\Psi)\subseteq\ell^2\setminus\ker(T_\Psi)$. This,
together with the fact that $S_\Psi$ is an invertible operator,
implies that $S_\Psi^{-1}T_{{\Psi}}((d_n)_n)\neq 0$. Now, if $f'$
is an element of $\h$ such that $(c_n)_n=U_\Psi f'$, then
\begin{eqnarray*}
T_{\widetilde{\Psi}}U_{\widetilde{\Phi}}T_\Phi U_\Psi(f')&=&
S_\Psi^{-1}T_{{\Psi}}U_{{\Phi}}S_\Phi^{-1}T_\Phi U_\Psi(f')\\
&=&S_\Psi^{-1}T_{{\Psi}}U_{{\Phi}}S_\Phi^{-1}T_\Phi((c_n)_n)\\
&=&S_\Psi^{-1}T_{{\Psi}}U_{{\Phi}}S_\Phi^{-1}T_\Phi(U_\Phi(g_1)+(d_n)_n)\\
&=&S_\Psi^{-1}T_{{\Psi}}(U_\Phi(g_1))+S_\Psi^{-1}T_{{\Psi}}U_{{\Phi}}S_\Phi^{-1}T_\Phi((d_n)_n)\\
&=&S_\Psi^{-1}T_{{\Psi}}((c_n)_n-(d_n)_n)+S_\Psi^{-1}T_{{\Psi}}U_{{\Phi}}S_\Phi^{-1}T_\Phi((d_n)_n)\\
&=&S_\Psi^{-1}T_{{\Psi}}((c_n)_n-(d_n)_n)\\
&=&f'-S_\Psi^{-1}T_{{\Psi}}((d_n)_n)\neq f',
\end{eqnarray*}
and this completes the proof.$\e$\\


\section{\large\bf Application to Gabor frames}

In order to state the results of this section we need to recall
the definition and some basic results on the duality conditions
for a pair of Gabor systems which have found more and more
applications in modern life, signal analysis and many other parts
of applied mathematics.

A {\it Gabor frame} is a frame for $L^2({\Bbb R})$ of the form
${\cal G}:=(E_{mb}T_{na}g)_{m,n\in{\Bbb Z}}$, where $a, b>0$ are
given, $g\in L^2({\Bbb R})$ is a fixed function,
$T_{na}f(x)=f(x-na)$ and $E_{mb}f(x)=e^{2imbx}f(x)$ for all $f\in
L^2({\Bbb R})$. In view of \cite[Theorem 9.1.12]{c}, the sequence
${\cal G}$  can only be a frame if $ab\leq 1$, but it is not a
sufficient condition. Moreover, if $\cal G$ is a frame and $ab<1$,
then there exists infinitely many $g^d$ in $L^2({\Bbb R})$ such
that we have the following reconstruction formula for each $f\in
L^2({\Bbb R})$
$$f=\sum_{m,n\in{\Bbb Z}}\big<f,E_{mb}T_{na}g^d\big>E_{mb}T_{na}g;$$
That is, the Gabor frames ${\cal G}$ and ${\cal
G}^d=(E_{mb}T_{na}g^d)_{m,n\in{\Bbb Z}}$ are dual frames. But the
standard choice of $g^d$ is $S_{\cal G}^{-1}g$, where $S_{\cal
G}:L^2({\Bbb R})\rightarrow L^2({\Bbb R})$ defined by
$$S_{\cal G}f=\sum_{m,n\in{\Bbb Z}}\big<f,E_{mb}T_{na}g\big>E_{mb}T_{na}g,
\quad\quad\quad (f\in L^2({\Bbb R}))$$ is the frame operator of
$\cal G$. There are some interesting results for duality
conditions of a Gabor frame $\cal G$, provided that the function
$g\in L^2({\Bbb R})$ has compact support, for more details see
\cite{dg1,c,dg2,a3}.

We begin the presentation of the results of this section by
recalling the following two results from \cite{dg1,dg2,jan}, each
of which will be required in our present investigation.

\begin{proposition} Let $N\in{\Bbb N}$ and $g\in L^2({\Bbb
R})$ be a function with support in $[0,N]$. Let also, $a, b>0$ be
given, $b\leq\frac{1}{N}$ and there exists $m, M>0$ such that
$$am\leq{\verb"G"}(x):=\sum_{n\in{\Bbb Z}}|g(x-na)|^2\leq bM\;\quad a.e.~x\in{\Bbb R}.$$
Then ${\cal G}=(E_{mb}T_{na}g)_{m,n\in{\Bbb Z}}$ is a frame for
$L^2({\Bbb R})$ with frame bounds $m, M$. Moreover, the frame
operator $S_{\cal G}$ and its inverse $S_{\cal G}^{-1}$ are given
by
\begin{eqnarray*}
S_{\cal G}(f)=\frac{b}{\verb"G"}f\quad\quad{\hbox{and}}\quad\quad
S_{\cal G}^{-1}(f)=\frac{\verb"G"}{b}f,
\end{eqnarray*}
for all $f\in L^2({\Bbb R})$.
\end{proposition}


The duality condition for a pair of Gabor systems ${\cal
G}=(E_{mb}T_{na}g)_{m,n\in{\Bbb Z}}$ and ${\cal
G}^d=(E_{mb}T_{na}g^d)_{m,n\in{\Bbb Z}}$ is presented by Janssen
\cite{jan} as follows.

\begin{lemma}
Two Bessel sequences ${\cal G}=(E_{mb}T_{na}g)_{m,n\in{\Bbb Z}}$
and ${\cal G}^d=(E_{mb}T_{na}g^d)_{m,n\in{\Bbb Z}}$ form dual
frames for $L^2({\Bbb R})$ if and only if
$$\sum_{k\in{\Bbb Z}}\overline{g(x-n/b-kn)}h(x-ka)=b\delta_{n,0}$$
for almost everywhere $x$ in $[0,a]$.
\end{lemma}


The proof of the next theorem can be found  in the paper by
Christensen and Kim \cite{dg2}.

\begin{theorem}\label{400}
Let $N\in{\Bbb N}$, $b\in(0,\frac{1}{2N-1}]$ and $g$ be a bounded
real-valued function for which
\begin{eqnarray*}
\sum_{n\in{\Bbb Z}}g(x-n)=1\quad\quad{\hbox{and}}\quad\quad {\hbox
{supp}}~g\subseteq[0,N].
\end{eqnarray*}
Then the function $g_1^d$ and $g_2^d$ defined by
\begin{eqnarray*}
g_1^d(x)=bg(x)+2b\sum_{n=1}^{N-1}g(x+n)\quad\quad {\hbox{and}}
\quad\quad g_2^d(x)=\sum_{n=-N+1}^{N-1}a_ng(x+n),
\end{eqnarray*}
where
\begin{eqnarray*}
a_0=b\quad\quad {\hbox{and}}\quad\quad a_n+a_{-n}=2b\quad
{\hbox{for~ each}}~~ n=1, 2,\cdots, N-1
\end{eqnarray*}
generate two dual frames ${\cal
G}_1^d=(E_{mb}T_{n}g_1^d)_{m,n\in{\Bbb Z}}$ and ${\cal
G}_2^d=(E_{mb}T_{n}g_2^d)_{m,n\in{\Bbb Z}}$ for ${\cal
G}=(E_{mb}T_{n}g)_{m,n\in{\Bbb Z}}$.\\
\end{theorem}


Here it should be noted that, Lemma 9.3.1 of \cite{c} guarantees
that there are infinitely many operator $\A$ on $L^2({\Bbb R})$
such that it commute with $E_{\pm b}$, $T_{\pm a}$ and
$\|Id_{L^2({\Bbb R})}-\A\|<1$, where $a$ and $b$ are positive real
numbers; In fact, each operator of the form
$$S_{\cal L}f=\sum_{m,n\in{\Bbb Z}}\big<f,E_{mb}T_{na}l\big>E_{mb}T_{na}l,
\quad\quad (f\in L^2({\Bbb R}))$$ gives an operator $\mathcal A$
on $L^2({\Bbb R})$ such that $\|Id_{L^2({\Bbb R})}-\A\|<1$ and it
commute with $E_{\pm b}$ and $T_{\pm a}$,
 where $l$ is a function in $L^2({\Bbb R})$ such that the sequence
 ${\cal L}=(E_{mb}T_{na}l)_{m,n\in{\Bbb Z}}$ is a frame for
 $L^2({\Bbb R})$. Moreover, we should mention that if an
 operator $\cal A$ commute with $E_{\pm b}$ and $T_{\pm
a}$, then ${\cal A}$ and its adjoint commute with $E_{mb}$ and
$T_{na}$ for all $m, n\in{\Bbb Z}$. Moreover, we have the
following example.


\begin{example}\label{ex1}
{\rm Assume that ${\cal G}=(E_{mb}T_{na}g)_{m,n\in{\Bbb Z}}$ and
${\cal G}^{ad}=(E_{mb}T_{na}g^{ad})_{m,n\in{\Bbb Z}}$ are
approximately dual frames. If we set ${\cal A}:=T_{{\cal
G}}U_{{\cal G}^{ad}}$, then $\|Id_{L^2({\Bbb R})}-{\cal A}\|<1$
and by a method similar to the proof of Lemma 9.3.1 of \cite{c}
one can easily see that ${\cal A}$ and its adjoint commute with
$E_{mb}$ and $T_{na}$ for all $m, n\in{\Bbb Z}$.}
\end{example}


Now, in view of Corollary \ref{200}, we have the following result
for the construction of approximately dual frames associated with
a given Gabor frame and given operator.

\begin{proposition}\label{300}
Let ${\cal G}=(E_{mb}T_{na}g)_{m,n\in{\Bbb Z}}$ be a frame for
$L^2({\Bbb R})$ and let ${\cal G}^d=(E_{mb}T_{na}g^d)_{m,n\in{\Bbb
Z}}$ be an arbitrary dual Gabor frame of $\cal G$. Suppose also
that $\cal A$ is an operator on $L^2({\Bbb R})$ such that
$\|Id_{L^2({\Bbb R})}-\A\|<1$ and it commutes with $E_{\pm b}$ and
$T_{\pm a}$. Then the family ${\cal
G}^{ad}=(E_{mb}T_{na}g^{ad})_{m,n\in{\Bbb Z}}$ is an approximately
dual frame of $\cal G$ for which $T_{\cal G}U_{{\cal G}^{ad}}=\A$,
where $g^{ad}=\A^* S_{\cal G}^{-1}g-g+S_{\cal G}(g^d)$.
\end{proposition}


It is well known from \cite{jan2} that in the case where $a\leq
c\leq 1$, then ${\cal
G}_{a,c}=(E_{m}T_{na}\chi_{[0,c)})_{m,n\in{\Bbb Z}}$ is a Gabor
frame, where $\chi_{[0,c)}$ denotes the characteristic function of
the interval $[0,c)$ on $\Bbb R$. Moreover, it was shown by
Hasankhani-Fard and Dehghan \cite[Corollary 2.1.]{dehghan2} that
two Bessel sequence ${\cal
G}_{a,c}=(E_{m}T_{na}\chi_{[0,c)})_{m,n\in{\Bbb Z}}$ and ${\cal
G}^d_{a,c'}=(E_{m}T_{na}\chi_{[0,c')})_{m,n\in{\Bbb Z}}$ are dual
frames for $L^2({\Bbb R})$ if and only if $c\leq 1$, $c'\leq 1$
and $a=\min\{c,c'\}$. From this, with the aid of Proposition
\ref{300} above, we can obtain the following explicit construction
of approximately dual frames associated with the Gabor frame
${\cal G}_{a,c}=(E_{m}T_{na}\chi_{[0,c)})_{m,n\in{\Bbb Z}}$ for
certain choice of $a$ and $c$. In the sequel, for Gabor frame
${\cal G}_{a,c}=(E_{m}T_{na}\chi_{[0,c)})_{m,n\in{\Bbb Z}}$, the
notation $S_{a,c}$ is used to denote its frame operator and
$\M_{a,c}$ denotes its Bessel bound.

\begin{example}
{\rm Let $a$, $c$, $c'$ and $c''$ be positive numbers such that
\begin{eqnarray*}
a\leq c, c'~{\hbox{and}}~ c''\leq
1,\quad\quad\quad{\hbox{and}}\quad\quad\quad a=\min\{c,c'\}.
\end{eqnarray*}

(a) Then two frame ${\cal
G}_{a,c}=(E_{m}T_{na}\chi_{[0,c)})_{m,n\in{\Bbb Z}}$ and ${\cal
G}_1^{ad}=(E_{m}T_{na}g_1^{ad})_{m,n\in{\Bbb Z}}$ are
approximately dual frames such that $T_{{\cal G}_{a,c}}U_{{\cal
G}^{ad}}=\frac{1}{\M_{a,c''}}S_{a,c''}$, where
$$g_1^{ad}=\frac{1}{\M_{a,c''}}\;S_{a,c''}\;
S_{a,c}^{-1}(\chi_{[0,c)})-\chi_{[0,c)}+S_{a,c}(\chi_{[0,c')}).$$

(b) Suppose that $\cal A$ is an operator on $L^2({\Bbb R})$ such
that $\|Id_{L^2({\Bbb R})}-\A\|<1$ and it commutes with $E_{\pm
1}$ and $T_{\pm a}$. Then the sequence ${\cal
G}_2^{ad}=(E_{m}T_{na}g_2^{ad})_{m,n\in{\Bbb Z}}$ is an
approximately dual frame of ${\cal
G}_{a,c}=(E_{m}T_{na}\chi_{[0,c)})_{m,n\in{\Bbb Z}}$ for which
$T_{{\cal G}_{a,c}}U_{{\cal G}^{ad}}=\A$, where $$g_2^{ad}=\A^*\;
S_{a,c}^{-1}(\chi_{[0,c)})-\chi_{[0,c)}+S_{a,c}(\chi_{[0,c')}).$$

(c) Let ${\cal G}=(E_{m}T_{na}g)_{m,n\in{\Bbb Z}}$ be a frame for
$L^2({\Bbb R})$ and let ${\cal G}^d=(E_{m}T_{na}g^d)_{m,n\in{\Bbb
Z}}$ be an arbitrary dual Gabor frame of $\cal G$. Suppose also
that
$$g_{a,e}^{ad}=\frac{1}{\M_{a,e}}S_{a,e}\;S_{\cal G}^{-1}g-g+S_{\cal G}(g^d),$$
then the sequence ${\cal
G}_{a,e}^{ad}=(E_{m}T_{na}g_{a,e}^{ad})_{m,n\in{\Bbb Z}}$ is an
approximately dual frame of $\cal G$ for which $T_{\cal G}U_{{\cal
G}_{a,e}^{ad}}=\frac{1}{\M_{a,e}}S_{a,e}$, where $e$ can be any positive real number such that
$a\leq e\leq 1$.\\ }
\end{example}


Recall from \cite[Section 6.1]{c} that for each $N\in{\Bbb N}$ the
$B$-splines $B_N$ are given inductively by
\begin{eqnarray*}
B_1=\chi_{[0,1]}\quad\quad{\hbox{and}}\quad\quad B_{N+1}=B_N\ast
B_1.
\end{eqnarray*}
The $B$-spline $B_N$ has support on the interval $[0,N]$.
Furthermore, for each $N\in{\Bbb N}$, the sequence ${\cal
G}_N=(E_{mb}T_{n}B_N)_{m,n\in{\Bbb Z}}$ is a Gabor frame with the
frame operator
\begin{equation}\label{310}
S_{N}f=\sum_{m,n\in {\Bbb
Z}}\big<f,E_{mb}T_{n}B_N\big>E_{mb}T_{n}B_N,\quad\quad (f\in
L^2({\Bbb R}))
\end{equation}
where $b\in(0,\frac{1}{2N-1}]$. Moreover, Theorem \ref{400} and
\cite[Theorem 6.1.1]{c} implies that the functions $B_{1,N}^d$ and
$B_{2,N}^d$ defined by
\begin{equation}\label{201}
B_{1,N}^d(x)=bB_N(x)+2b\sum_{n=1}^{N-1}B_N(x+n)
\end{equation}
and
\begin{equation}\label{202}
B_{2,N}^d(x)=\sum_{n=-N+1}^{N-1}a_nB_N(x+n),
\end{equation}
where
\begin{eqnarray*}
a_0=b\quad\quad {\hbox{and}}\quad\quad
a_n+a_{-n}=2b\quad {\hbox{for~ each}}~ n=1, 2,\cdots, N-1,
\end{eqnarray*}
generate two dual frames ${\cal
G}_{1,N}^d=(E_{mb}T_{n}B_{1,N}^d)_{m,n\in{\Bbb Z}}$ and ${\cal
G}_{2,N}^d=(E_{mb}T_{n}B_{2,N}^d)_{m,n\in{\Bbb Z}}$ for ${\cal
G}_N=(E_{mb}T_{n}B_N)_{m,n\in{\Bbb Z}}$.\\


As an application of Proposition \ref{300} with the aid of Example
6.1 of \cite{app} one can easily obtain the following result for
the second order B-spline.

\begin{example}
{\rm Let $\varphi(x)=e^{-4x^2}$ and
$$g(x)=\frac{15.1}{315}\;\frac{1}{\sum_{n\in{\Bbb Z}}|B_8(2.36(x+n))|^2}\;B_8(2.36x).$$
Then, by Example 6.1 of \cite{app}, we deduce that the frames
${\Phi}=(E_{0.1m}T_{n}\varphi)_{m,n\in{\Bbb Z}}$ and ${\cal
G}=(E_{0.1m}T_{n}g)_{m,n\in{\Bbb Z}}$ are approximately dual
frames and $\|Id_{L^2({\Bbb R})}-{\cal A}\|\leq 0.009$. If now we
consider ${\cal A}:=T_{\Phi}U_{{\cal G}}$, $$B_2^{ad}=\A^*
S_2^{-1}B_2-B_2+S_2(B_{1,2}^d)$$ and
$$B_{1,2}^d(x)=0.1B_2(x)+0.2B_2(x+1),$$ then the sequence
${\cal G}_2^{ad}=(E_{0.1m}T_{n}B_2^{ad})_{m,n\in{\Bbb Z}}$ is an
approximately dual frame of ${\cal G}_2$ for which $T_{{\cal
G}_2}U_{{\cal G}_2^{ad}}=\A$. In particular, the
approximation rate is $\varepsilon=0.009$. }
\end{example}


Finally, with the above notations on the N-order B-spline, we have the
following practical example.

\begin{example}
{\rm  Let $N\in{\Bbb N}$, $b\in (0,\frac{1}{2N-1}]$ and $\cal A$
be an operator on $L^2({\Bbb R})$ such that $\|Id_{L^2({\Bbb
R})}-\A\|<1$ and it commutes with $E_{\pm b}$ and $T_{\pm 1}$. Let
also $S_{N}$ be given by equation {\rm(\ref{310})} and $\M_N$ be
the Bessel bound of the frame ${\cal
G}_N=(E_{mb}T_{n}B_N)_{m,n\in{\Bbb Z}}$.

(a) Then the function $B_{1,N}^{ad}$ and $B_{2,N}^{ad}$ defined by
$$
B_{1,N}^{ad}(x)=\A^* S_{N}^{-1}B_N-B_N+S_{N}(B_{1,N}^d)
$$
and
$$
B_{2,N}^{ad}(x)=\A^* S_{N}^{-1}B_N-B_N+S_{N}(B_{2,N}^d),
$$
generate two approximately dual frames ${\cal
G}_{1,N}^{ad}=(E_{mb}T_{n}B_{1,N}^{ad})_{m,n\in{\Bbb Z}}$ and
${\cal G}_{2,N}^{ad}=(E_{mb}T_{n}B_{2,N}^{ad})_{m,n\in{\Bbb Z}}$
for ${\cal G}_N=(E_{mb}T_{n}B_N)_{m,n\in{\Bbb Z}}$ such that
$T_{{\cal G}_N}U_{{\cal G}_{1,N}^{ad}}=\A=T_{{\cal G}_N}U_{{\cal
G}_{2,N}^{ad}}$, where $B_{1,N}^d$ and $B_{2,N}^d$ are given by
equations {\rm(\ref{201})} and {\rm(\ref{202})}.

(b) Let ${\cal G}=(E_{mb}T_{n}g)_{m,n\in{\Bbb Z}}$ be a frame for
$L^2({\Bbb R})$ and let ${\cal G}^d=(E_{mb}T_{n}g^d)_{m,n\in{\Bbb
Z}}$ be an arbitrary dual Gabor frame of $\cal G$. Suppose that
$N$ is an arbitrary element of ${\Bbb N}$ and
$$g_N^{ad}=\frac{1}{\M_N}S_{N}\; S_{\cal G}^{-1}g-g+S_{\cal G}(g^d).$$
Then the family ${\cal
G}_N^{ad}=(E_{mb}T_{n}g_N^{ad})_{m,n\in{\Bbb Z}}$ is an
approximately dual frame of $\cal G$ for which $T_{\cal G}U_{{\cal
G}_N^{ad}}=\frac{1}{\M_N}S_{N}$.\\ }
\end{example}


\vspace{3mm}

\noindent {\sc Hossein Javanshiri}\\
Department of Mathematics,
Yazd University,
P.O. Box: 89195-741, Yazd, Iran\\
E-mail: h.javanshiri@yazd.ac.ir\\

\end{document}